\documentclass[11pt]{article}
\usepackage{lipsum}
\usepackage{amsfonts}
\usepackage{graphicx}
\usepackage{epstopdf}
\usepackage{amsmath}
\usepackage{amsfonts}

\usepackage{moreverb}
\usepackage{stmaryrd}
\usepackage{float}
\usepackage{subfig}
\usepackage{algorithm}
\usepackage{algorithmic}
\usepackage{bbm}
\usepackage{nicefrac}
\usepackage{hyperref}

\usepackage{graphicx}
\usepackage{epstopdf}

\newcommand{\bs}[1]{\ensuremath{\boldsymbol{#1}}}
\newcommand{\bb}[1]{\ensuremath{\mathbb{#1}}}

\newcommand{\trans}{\mathrm{T}}
\newcommand{\dd}{\mathrm{d}}

\title{Efficient Explicit Time Stepping of High Order Discontinuous Galerkin Schemes for Waves}
\author{S.~Schoeder (1), K.~Kormann (2), W.A.~Wall (1), M.~Kronbichler (1)\\
\small{(1) Institute for Computational Mechanics, Technical University of Munich,}\\ \small{Garching, Germany}\\ \small{(2) Max-Planck-Institut f\"ur Plasmaphysik, Garching, Germany}}
\date{\today}

\begin{document}

\maketitle

\section*{Abstract}
This work presents algorithms for the efficient implementation of discontinuous Galerkin methods with explicit time stepping for acoustic wave propagation on unstructured meshes of quadrilaterals or hexahedra. A crucial step towards efficiency is to evaluate operators in a matrix-free way with sum-factorization kernels. The method allows for general curved geometries and variable coefficients. 
Temporal discretization is carried out by low-storage explicit Runge--Kutta schemes and the arbitrary derivative (ADER) method. For ADER, we propose a flexible basis change approach that combines cheap face integrals with cell evaluation using collocated nodes and quadrature points. Additionally, a degree reduction for the optimized cell evaluation is presented to decrease the computational cost when evaluating higher order spatial derivatives as required in ADER time stepping.
We analyze and compare the performance of state-of-the-art Runge--Kutta schemes and ADER time stepping with the proposed optimizations. ADER involves fewer operations and additionally reaches higher throughput by higher arithmetic intensities and hence decreases the required computational time significantly.  Comparison of Runge--Kutta and ADER at their respective CFL stability limit renders ADER especially beneficial for higher orders when the Butcher barrier implies an overproportional amount of stages. 
Moreover, vector updates in explicit Runge--Kutta schemes are shown to  take a substantial amount of the computational time due to their memory intensity. 

\noindent
\textbf{Keywords:} 
  matrix-free methods,
  sum factorization,
  discontinuous Galerkin methods,
  arbitrary high-order,
  explicit time stepping,
  high performance computing

\noindent
\textbf{AMS:}
  65M60, % Finite elements, Rayleigh-Ritz and Galerkin methods, finite methods
  65Y20, % Complexity and performance of numerical algorithms
  68Q25, % Analysis of algorithms and problem complexity
  68W40  % Analysis of algorithms

\section{Introduction}

Discontinuous Galerkin (DG) methods are highly attractive for solving partial differential equations of transport character because they combine the robustness of finite volume methods with the accuracy of finite element methods on unstructured meshes~\cite{cks00,hw08}. 
An essential  component is the implementation of DG on modern hardware, which is the focus of this work. Modern processors favor arithmetically intense algorithms such as dense matrix multiplication over patterns typical for many solvers of partial differential equations that include loads of a streaming character with much lower arithmetics. 
Thus, algorithmic complexities alone are not enough to judge a method's efficiency as the achievable performance of e.g.~a sparse matrix-vector product in terms of floating point operations per second can be almost two orders of magnitude below the advertised peak performance \cite{roofline}.

This paper derives and analyzes two time stepping schemes for the discontinuous Galerkin (DG) method for the acoustic wave equation that are significantly more efficient than traditional matrix-based algorithms. 
The first scheme is based on explicit low-storage Runge--Kutta methods, and the second one on arbitrary derivative (ADER) time stepping that relies on a temporal Taylor expansion of the solution and expresses temporal derivatives by spatial derivatives using the Cauchy--Kowalevski (also called Lax--Wendroff) procedure. ADER has been successfully applied in the context of finite volume and DG methods~\cite{dk06,dkt07,sdm04}. An advantage of ADER over explicit Runge--Kutta schemes is that ADER is not restricted by the Butcher barriers and convergence orders beyond four are not overproportionally expensive.
Recently, ADER was combined with the hybridizable discontinuous Galerkin (HDG) method~\cite{skw17}. HDG had originally been introduced in the context of implicit time integration \cite{cgl09}, but it was later found that it is highly attractive in combination with explicit time integration outperforming implicit time integration by orders of magnitude \cite{ksmw15}. A characteristic property of HDG discretizations is that a superconvergent solution can be constructed at very low additional cost. In~\cite{skw17}, a reconstruction procedure is presented that allows to combine ADER and HDG while maintaining the superconvergence property.

Previous work on ADER-DG (as in~\cite{bhrb14}) relies on triangles or tetrahedra assuming constant coefficients and straight lined boundaries. The operator evaluation in~\cite{bhrb14} is carried out based on an element matrix with a theoretical complexity per degree of freedom of $\mathcal O(k^d)$ in the degree $k$ for all spatial dimensions $d$. Here, we propose an ADER-DG formulation for quadrilaterals and hexahedra for variable coefficients and curved geometries. We  use a matrix-free operator evaluation relying on fast quadrature with sum-factorization kernels with a theoretical complexity per degree of freedom of $\mathcal O(dk)$. The techniques of sum factorization with fast quadrature have been established by the spectral element community~\cite{Deville02,ks05,kopriva,Orszag80} but are also popular in the DG community for explicit time integration~\cite{Hindenlang12}. Advances in computer architecture have rendered the matrix-free evaluation, originally targeting high orders beyond around five, also highly competitive at moderate orders, outperforming the memory bandwidth-limited sparse matrix-vector product for second and higher degree polynomials~\cite{brown10,kk12}. 
In terms of algorithmic layout, the sole reduction of arithmetic operations is not advantageous if the memory bandwidth is the performance limiting factor. In this work, we show that an operator evaluation with sum factorization as in explicit Runge--Kutta schemes is memory bandwidth bound, despite its clear improvement over matrix-based operator evaluation. ADER replaces the global operator application in each Runge--Kutta stage by one global operator application and a completely element-local evaluation routine, the Taylor--Cauchy--Kowalevski procedure, which allows to perform more computations on data read from the global solution vectors. We show that ADER does not only employ fewer operations but also supplies a higher arithmetic intensity which is beneficial on modern cache-based hardware. 

Another aspect significant for performance is the choice of the shape function nodal points and the choice of the quadrature rule. In case nodal points and quadrature points coincide, interpolation of the solution to quadrature points is avoided and computational expense is decreased. This approach is well known in the context of spectral elements. 
Usage of the Gauss-Lobatto points for the definition of the nodal points and for integration was shown to degrade the accuracy of the mass matrix and its inverse~\cite{Durufle09,teukolsky15}, though. Consistent Gaussian quadrature instead yields full accuracy. A drawback of nodes in the Gauss points, however, is that the flux evaluation on element faces requires an extrapolation accessing all degree of freedom values of both adjacent elements because there are no node points on the faces. For the flux evaluation, nodal points on the element faces ensure that only $(k+1)^{d-1}$ instead of $(k+1)^d$ values must be accessed. 
We propose a new algorithmic method that changes the polynomial basis and its nodes on the fly depending on the quantity to be evaluated. The standard DG global derivative operator including flux evaluations relies on a Lagrange basis with Gauss-Lobatto points while the ADER specific element local Taylor--Cauchy--Kowalevski procedure relies on a Lagrange basis with nodes in collocated Gauss points. Thereby, we are able to combine cheap element evaluation and flux evaluation with minimal data access.
Despite this optimization concerning node and quadrature choice, we propose a second optimization concerning the efficient evaluation of higher order spatial derivatives required in the Taylor--Cauchy--Kowalevski procedure. Calculation of first order spatial derivatives and successive projection to a lower order basis in combination with the collocated node and quadrature points minimize computational work.

This article is structured as follows. The physical problem and its spatial and temporal discretization are given in Section~\ref{sec:basics}. 
The basic algorithmic building blocks are described in the first part of Section~\ref{sec:algo}. 
A quantitative study on the throughput for basis with or without nodes on element surfaces is given in~\ref{sec:prelim}, which motivates the basis switching approach presented in~\ref{sec:basis-change}. The optimization relying on reduced polynomial spaces for higher spatial derivatives is shown in Section~\ref{sec:basis-reduction}.
A detailed performance analysis in terms of theoretically derived operation counts, throughput, the roofline model, computational timings, and scalability is given in Section~\ref{sec:results}. We conclude in Section~\ref{sec:concl}.

\section{The Problem and its Discretization}\label{sec:basics}
The acoustic wave equation in terms of the solution vector function $\bs{u}=[\bs{v},p]$ summarizing the particle velocity $\bs{v}$ and the pressure deviation $p$ is written as
\begin{equation}
\frac{\partial}{\partial t} \bs{u} +\bs{S} \bs{u} = \bs{0} \label{eq:general-hyp-pde}
\end{equation}
with the matrix $\bs{S}$ containing parameters and first order spatial derivatives in the form 
\begin{align*}
\bs{S} = \left[ \begin{array}{cc} 0 & \frac{1}{\rho} \nabla \\ c^2\rho \nabla \cdot & 0 \end{array} \right]
\end{align*}
with speed of sound $c$ and mass density $\rho$. Most algorithmic components developed in this work can also be applied to other hyperbolic operators $\bs{S}$.
The solution on a $d$-dimensional domain $\Omega\in\mathbb{R}^d$ is sought with respect to initial conditions $\bs{u}(t=0)=\bs{u}_0$ and boundary conditions. To find an approximate solution to problem~\eqref{eq:general-hyp-pde}, first the spatial discretization and subsequently the temporal discretization are applied.

\subsection{The Discontinuous Galerkin Method}\label{sec:dg}
The domain is tesselated into a triangulation $\mathcal{T}$ of quadrilaterals or hexahedra. For discretization with DG, equation~\eqref{eq:general-hyp-pde} is multiplied with test functions $\bs{w}$, integrated over all elements of $\mathcal{T}$ and integrated by parts. The resulting problem reads: Find $\bs{u}$ in the approximation space spanned by element-wise polynomials of tensor degree $k$, constructed as the tensor product of one-dimensional functions, such that
\begin{align}
\left( \bs{w}, \frac{\partial}{\partial t} \bs{u} \right)_\mathcal{T} - \left( \bs{S}^\trans\bs{w}, \bs{u} \right)_\mathcal{T} + \left\langle \bs{w}, \hat{\bs u} \right\rangle_{\partial \mathcal{T}} = 0, \label{eq:weak-form}
\end{align}
with $(\cdot,\cdot)_\mathcal{T}$ and $\langle\cdot,\cdot\rangle_{\partial \mathcal{T}}$ denoting the $L_2$ inner product over the $d$ and $(d-1)$ dimensional domains $\mathcal{T}$ and ${\partial \mathcal{T}}$, respectively.
The test functions are from the same function space as the solution functions. The boundary of the tesselation $\partial \mathcal{T}$ comprises all element boundaries both at the interior faces and the boundary of the physical domain $\Omega$. The variable $\hat {\bs u}$ is the numerical flux for evaluation of the element surface integrals by relating quantities from both adjacent elements. Several possibilities for the flux definition exist, such as central fluxes, upwinding, HDG fluxes, and many more. The choice depends on the specific problem and intended numerical properties.

Expressing the solution $\bs{u}$ as product of ansatz functions summarized in a matrix $\bb{N}$ and values of the degrees of freedom $\bs{U}$, i.e.,
\begin{align*}
\bs{u}=\bb{N}\bs{U},
\end{align*}
 the time continuous space discrete system is derived
\begin{equation}
\bb{M} \frac{\partial}{\partial t} \bs{U} + \bb{K} \bs{U} = \bs{0}, \label{eq:matrix-tc-sd}
\end{equation}
with a mass matrix $\bb{M}$ and a stiffness matrix $\bb{K}$. The stiffness matrix entries depend on the choice of the flux definition $\hat{\bs u}$. 

\subsection{Explicit Runge--Kutta Time Integration}\label{sec:rk}
Time discretization of equation~\eqref{eq:matrix-tc-sd} using a general explicit Runge--Kutta method involves summation over the $s$ Runge--Kutta stages and repeated operator application
\begin{equation}
\bs{U}_{t_{i+1}} = \bs{U}_{t_{i}} + \Delta t \sum_{j=1}^s b_j k_j \quad \text{ and } \quad k_j = -\bb{M}^{-1}\bb{K} \left(\bs{U}_{t_{i}} + \Delta t \sum_{l=1}^{j-1} a_{jl} k_l \right),
\end{equation}
with the coefficients $a_{jl}$ and $b_j$ given in Butcher tableaus for specific Runge--Kutta schemes. A well known explicit Runge--Kutta scheme is the classical Runge--Kutta of order four with four stages. Optimized variants, e.g. low-storage schemes, which were originally introduced to reduce the memory consumption but that we use primarily for their lower the memory transfer, are given in~\cite{kcl00}.

In this work, the low-storage schemes with two registers of order four with five stages LSRK4(5) and of order five with nine stages LSRK5(9) from~\cite{kcl00} will be used. 
According to the Butcher barriers, there is no five stage Runge--Kutta scheme of order five. The minimum number of stages is six for accuracy order five. With the additional conditions imposed by the requirement to obtain a low-storage scheme with two registers, the minimal number of stages is nine. An eight stage solution is presumed to be impossible~\cite{kcl00}. For an alternative variant of low-storage schemes, which was explicitly developed for DG, see e.g.~\cite{Desmet12}; these schemes are not considered in this work because their efficiency is not significantly better than for the schemes from~\cite{kcl00}.

\subsection{Arbitrary Derivative Time Integration}\label{sec:ader}
ADER explicit time integration has been presented in the context of DG spatial discretization and several hyperbolic problems, see~\cite{dk06,dprz16,sdm04}. Time discretization of the time-continuous system~\eqref{eq:matrix-tc-sd} with ADER relies on the idea to express temporal derivatives with spatial derivatives using the partial differential equation~\eqref{eq:general-hyp-pde} according to the Cauchy--Kowalevski procedure
\begin{align*}
\frac{\partial}{\partial t} \bs{u} = -\bs{Su} \qquad \Rightarrow \qquad \frac{\partial^j}{\partial t^j} \bs{u} = (-1)^j \bs{S}^j \bs{u}.
\end{align*}
The solution field is expanded into a truncated Taylor series in time around the last known state at $t=t_i$ and the above expression is inserted. A projection onto the degree of freedom values yields
\begin{align}
\bs{U} = \bb{M}^{-1} \sum_{j=0}^{k+1} \frac{(t-t_i)^j}{j!} (-1)^j\int_K\bb{N}^\trans \bs{S}^j\bb{N}\dd K \bs{U}_{t_i}, \label{eq:ader-U}
\end{align}
where the subscript $(\cdot)_{t_i}$ means evaluation at $t_i$. Last, equation~\eqref{eq:matrix-tc-sd} is integrated in time and equation~\eqref{eq:ader-U} is used to replace the time integral. The final scheme reads
\begin{align}
\bs{U}_{t_{i+1}} = \bs{U}_{t_i} - \bb{M}^{-1} \bb{K} \bb{M}^{-1} \sum_{j=0}^{k+1} \frac{(t_{i+1}-t_i)^{j+1}}{(j+1)!} (-1)^j\int_K\bb{N}^\trans \bs{S}^j\bb{N}\dd K \bs{U}_{t_i}. \label{eq:ader}
\end{align}
Just as every explicit time integration scheme, ADER time stepping is subject to time step restrictions according to the Courant--Friedrichs--Lewy (CFL) condition~\cite{cfl28} with the Courant number $\textsl{Cr}$
\begin{align*}
\textsl{Cr} = \frac{c \Delta t k^{1.5}}{h}
\end{align*}
in terms of the signal speed $c$, the time step size $\Delta t$, the polynomial degree of the shape functions $k$, and a characteristic element size $h$. The power $1.5$ on the polynomial degree is used to correct for higher order approximations as suggested in~\cite{ks05}. The stability restriction for ADER is generally stricter than for explicit Runge--Kutta methods, as mentioned in~\cite{dk06} where $50\%$ of the Runge--Kutta Courant number is given as rule of thumb. Time step stability is investigated in Section~\ref{sec:time-vs-acc}.

\subsection{The Hybridizable Discontinuous Galerkin Method}\label{sec:hdg}
The hybridizable discontinuous Galerkin method (HDG) in combination with explicit time stepping can be understood as a special choice of the flux $\hat{\bs u}$ in equation~\eqref{eq:weak-form}. For the acoustic wave equation, the HDG flux is
\begin{equation*}
\hat{\bs{u}}_\text{HDG}=\left[ \begin{array}{cc}
\frac{1}{2\rho} \left(p^+_h+p^-_h\right)\bs{n}^- & + \frac{1}{2\rho\tau} \left(\bs{n}^-\cdot\bs{v}^-_h -\bs{n}^-\cdot\bs{v}^+_h\right) \bs{n}^- \\
\frac{c^2\rho}{2}\left(\bs{n}^-\cdot\bs{v}^-_h - \bs{n}^+\cdot\bs{v}^+_h\right) &+\frac{c^2\rho\tau}{2} \left(p^-_h-p^+_h\right)
\end{array}\right],
\end{equation*}
with the superscripts $(\cdot)^+$ and $(\cdot)^-$ indicating evaluation of the two adjacent elements and $\bs{n}$ being the outward pointing normal vector. HDG introduces a stabilization parameter $\tau$ typically chosen as $\tau = \nicefrac{1}{c\rho}$. If a time integration scheme of sufficient accuracy is utilized, HDG offers the possibility to obtain superconvergent pressure results of order $k+2$ by means of a simple postprocessing step. To maintain the superconvergence property of HDG in combination with ADER time integration where spatial and temporal discretization are strongly interlinked, the method shown in equation~\eqref{eq:ader} is enhanced to
\begin{align}
\begin{split}
\bs{U}_{t_{i+1}} &= \bs{U}_{t_i} - (t_{i+1}-t_i) \bb{M}^{-1}\bb{K} \bs{U}_{t_i}  \\
&- \bb{M}^{-1} \bb{K} \bb{M}^{-1} \left( \sum_{j=1}^{k+1} \frac{(t_{i+1}-t_i)^{j+1}}{(j+1)!} (-1)^j\int_K\bb{N}^\trans \bs{S}^{j-1}\bb{N}\dd K \right) \bb{M}^{-1} \bb{K} \bs{U}_{t_i},
\end{split}\label{eq:hdg-ader}
\end{align}
replacing one element-local derivative $\bs{S}$ by a discrete DG derivative $\bb{M}^{-1}\bb{K}$. This procedure recovers a superconvergent scheme, see~\cite{skw17} for details.

\section{Algorithmic Developments}\label{sec:algo}

In the beginning of this section, we describe the basic underlying algorithmics on which our new developments build. Subsequently, we give a brief preliminary performance evaluation for different polynomial bases in Section~\ref{sec:prelim}. Thereby, we motivate a basis change as described in Section~\ref{sec:basis-change}, which is optimized in terms of a degree reduction in Section~\ref{sec:basis-reduction}. 

The evaluation of the integrals in the weak forms in equations~\eqref{eq:ader} or~\eqref{eq:hdg-ader} is performed by fast integration relying on sum factorization utilizing the structure of tensor product shape functions that are integrated with a tensor product quadrature rule. Throughout this work, we choose a Gaussian quadrature with $k+1$ points per direction for polynomials of degree $k$, which is enough to integrate bilinear forms with element-wise constant coefficients on affine element shapes exactly. On curved geometries, there is a possible integration error that is often subsumed in the errors from variational crimes \cite{brenner08}. In particular, our choice avoids the accuracy penalty of inexact Gauss--Lobatto quadrature on $k+1$ points as highlighted in \cite{Durufle09}. Note that more general hyperbolic problems with nonlinear terms can easily be integrated with more points to avoid aliasing effects in this setup.

For a function described by a basis of degree $k$ and $(k+1)^d$ coefficients, 
the interpolation onto $(k+1)^d$ quadrature points takes $(k+1)^{2d}$ additions and multiplications in a naive implementation without sum factorization.
Also, the evaluation of each component of the gradient takes $2(k+1)^{2d}$ operations. With sum factorization however, the work for the interpolation is reduced to $2d(k+1)^{d+1}$ operations.
The reduction of operations results from $d$ applications of one-dimensional interpolations of cost $2(k+1)^2$ that go through all $(k+1)^{d-1}$ lines of basis functions and quadrature points, respectively. In the remainder of this work, we refer to one application of a one-dimensional interpolation over all points, involving $(k+1)^{d+1}$ operations, as one ``tensor product kernel''. Note that these kernels can be cast as small matrix-matrix multiplications. For details, we refer to~\cite{kk17}.

In case of quadrature over Lagrange polynomials with nodes in the points of the quadrature formula (a so-called collocated setup \cite{kopriva}) the interpolation from coefficients to values in quadrature points is the identity operation and can be skipped. Thus, the evaluation of all $d$ components of the gradient only involves $d$ tensor product kernels for each of the partial derivatives, as opposed to $d^2$ tensor product kernels for the basic evaluation scheme.  
As will be elaborated in this work, the optimization of using a collocated bases may be premature as other factors, such as the cost of face integrals, may control the decision about which basis to choose. Similar kernels are also used for the summation over quadrature points when multiplying with test functions or test function gradients, see e.g.~\cite{kk17} for details.

For our experiments, we use a state-of-the-art implementation of sum factorization based on the \texttt{deal.II} finite element library with support for massively parallel computations and adaptively refined meshes with hanging nodes \cite{dealII85}. Since integration involves a series of heavy arithmetics, the use of vectorization (SIMD) is fundamental for getting optimal performance on current architectures. Following the concepts described in \cite{kk17,kronbichler17isc}, this work applies vectorization over several elements which was found to provide best performance on polynomial degrees up to at least 14 and reaches more than 50\% of the arithmetic peak performance when considered in isolation, which is an extremely high value for a code compiled from generic C++ code that contains $k$ as a (template) parameter that lets the compiler decide on the loop unrolling and register allocation.

\subsection{Efficient Face Integral Evaluation}\label{sec:prelim}

In this section, the effect of the polynomial basis functions on the evaluation of the derivative operator $\bb{K}$ and the inverse mass matrix $\bb{M}^{-1}$ are studied. Two contradictory factors are considered: The inverse mass matrix is efficiently evaluated if quadrature points and nodes of the polynomial basis coincide, which results in a diagonal mass matrix and hence simple inversion. On the other hand, the  derivative operator $\bb{K}$ includes the evaluation of element  but also face integrals. For the face integrals, data from both adjacent elements are required. For nodal polynomials with nodes on the element faces, only the values associated to the nodes on the faces must be accessed. If the nodes are only in the interior or polynomials are not nodal, all vector entries of both adjacent elements are required and an extrapolation to the faces must be carried out. 

\begin{figure}
\centering
\includegraphics[scale=1.0]{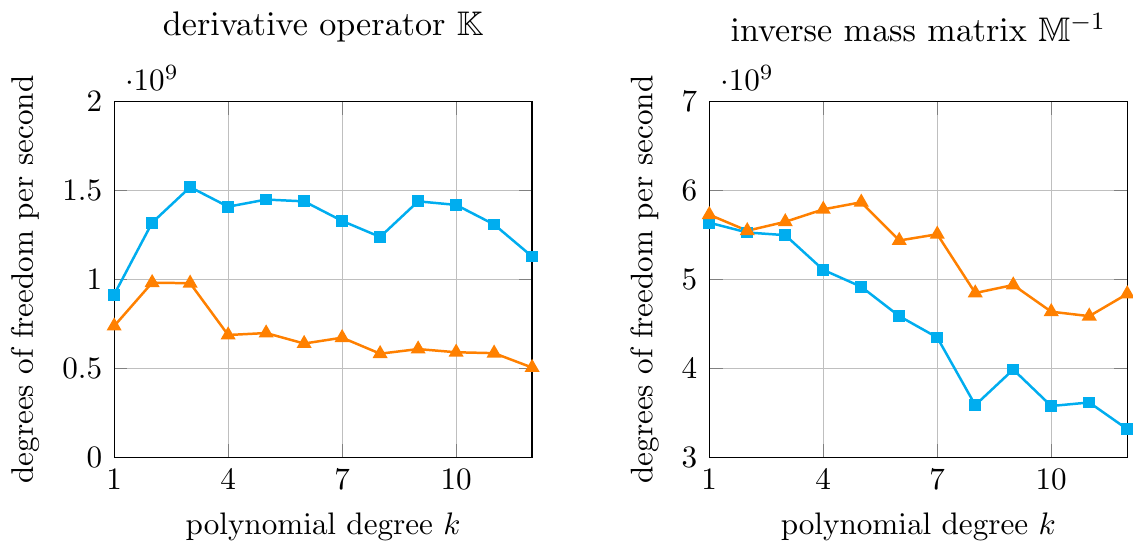}
\includegraphics[scale=1.0]{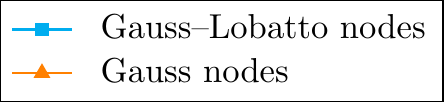}
\caption{Degrees of freedom processed per second for the derivative operator $\bb{K}$ and the inverse mass matrix $\bb{M}^{-1}$ on a Cartesian grid for polynomial degrees in $k\in\{1,\ldots,12\}$. Gaussian quadrature on $(k+1)^d$ points is used. In one case, the polynomial basis consists of Lagrange polynomials with nodes in Gauss--Lobatto points. In the other case, the Lagrange polynomials have their nodes in the Gauss points yielding a diagonal mass matrix.}
\label{fig:rk4-basis}
\end{figure}

Two variants are examined to demonstrate the effects of the aforementioned contradictory factors. Quadrature is always based on $(k+1)^d$ Gaussian points. In one case, a polynomial basis of Lagrange functions with nodes in the same Gauss points is used, in the other case a polynomial basis of Lagrange functions with nodes in Gauss--Lobatto points is used.  The first case results in a diagonal mass matrix, the second case satisfies the requirement to have nodes on the element faces. Note that the usage of Gauss--Lobatto points as nodes and quadrature points is not considered because it degrades accuracy~\cite{Durufle09,teukolsky15}. Fig.~\ref{fig:rk4-basis} plots the results in terms of degrees of freedom processed per second for a three-dimensional geometry consisting of $80^3$ Cartesian elements for polynomial degree $k\in\{1,2,3\}$ and $40^3$ Cartesian elements for $k\in\{4,\ldots,12\}$. Computations are run on the system specified in Table~\ref{tab:broadwell}.

The evaluation of the inverse mass matrix reaches a considerably higher throughput compared to the derivative operator $\bb{K}$ because it is completely element-local whereas $\bb{K}$ requires the evaluation of face integrals and therefore additional data access from neighboring elements. In fact, the throughput of the former is mostly limited by the memory bandwidth of loading the vector and writing back into the same vector, which is $6.2\cdot 10^{9}$ degrees of freedom per second for the system's memory throughput of 115 GB/s when counting $2\times 8(d+1)$ byte per polynomial (read and write) and $8$ bytes for access to the precomputed inverse entry of the diagonal, which is the same for all $(d+1)$ components. Following \cite{ksmw15}, the action of the inverse can be cheaply evaluated by sum-factorization kernels that first transform into an orthogonal basis (i.e., the Lagrange polynomials in the $(k+1)^d$ points of Gauss quadrature), apply the inverse diagonal mass matrix, and transform back to the original basis.
The comparison to the throughput for collocated nodes of shape functions and quadrature points shows that indeed the evaluation of the non-diagonal mass matrix operation is only up to 15\% slower compared to the diagonal mass matrix. Thus, the substantial arithmetic operations can be almost completely hidden behind the memory transfer. 

The face evaluation for the Lagrange basis in Gauss points does not only require to read the values for the degrees of freedom located on the face but all values of the adjacent elements to allow interpolation onto the face (or to write the face data into a separate global array with additional memory transfer as e.g.~done in \cite{hw08,Hindenlang12}). Comparing the throughput for Lagrange polynomials in Gauss--Lobatto points and Gauss points reveals a significant drop in performance for the derivative operator due to the increased memory access, including substantial indirect addressing components as described in \cite{kk17}. The direct comparison highlights that using polynomial bases with support points on the element faces for the global derivative operator is highly beneficial. This finding applies particularly to Runge--Kutta type time integrators, and it also holds for more general nonlinear operators $\bs{S}$.

We want to conclude this section by summarizing the main findings:
\begin{itemize}
\item The evaluation of the inverse mass matrix is memory bandwidth bound (especially for moderate order) and a change of basis is for free.
\item The throughput for the derivative operator is much lower than for the mass matrix, because values from both adjacent elements must be read in contrast to complete element-local evaluations. Combined with the fact that also the Jacobian of the mapping must be accessed, the memory bandwidth is reached earlier.
\item This effect is counteracted by usage of a polynomial basis with nodes on the element boundary, which will be further investigated in the following section.
\item Usage of collocated node and quadrature points reduces the number of operations.
\end{itemize}
From these conclusions, ADER is self-evidently motivated: Applications of the global derivative operator $\bb{K}$ are traded for element local evaluations in the Taylor--Cauchy--Kowalevski procedure.

\subsection{Flexible Basis Change}\label{sec:basis-change}
In the previous section, the effect of the memory bandwidth bound became apparent. Despite the memory access, the number of operations is a central quantity of interest to judge code efficiency. 

As shown in Fig.~\ref{fig:rk4-basis}, a set of shape functions where only $(k+1)^{d-1}$ functions evaluate to non-zero on each of the element faces is beneficial to reduce the vector access for face evaluation. For element local evaluations, however, this involves additional $d$ tensor product kernels per component to interpolate from the solution values to quadrature points as compared to the collocated node and quadrature points commonly used in spectral element solvers~\cite{kopriva}. 
Hence, different bases should be used for the different phases. For element local evaluations, a collocated basis with nodes in Gauss points (denoted by G) is the best choice while a basis with nodes in Gauss--Lobatto points (denoted by GL) is the best choice for the evaluation of face integrals. We propose to switch the basis for the different phases on the fly while looping over the elements for the Taylor--Cauchy--Kowaleski evaluation. 
The solution approximation is expressed either as
\begin{align*}
\bs{u}_\text{G} = \bb{N}_\text{G} \bs{U}_\text{G} \qquad \text{ or } \qquad \bs{u}_\text{GL} = \bb{N}_\text{GL} \bs{U}_\text{GL},
\end{align*}
with the matrices $\bb{N}_\text{G},\bb{N}_\text{GL}$ containing the shape functions evaluated at the respective nodes and the vectors $\bs{U}_\text{G} ,\bs{U}_\text{GL}$ containing the degree of freedom values for Gauss and Gauss--Lobatto points, respectively. 
Theoretically, there are no restrictions objecting to change the basis from one evaluation to the other. Since both spaces are of the same degree, the equality $\bs{u}_\text{G}=\bs{u}_\text{GL}$ holds in all cases.

For ADER, the collocated G basis is not only interesting for the application of the inverse mass matrix $\bb{M}^{-1}$ but especially for the Taylor--Cauchy--Kowalevski term
\begin{equation*}
\sum_{j=0}^{k+1} \frac{(t_{i+1}-t_i)^{j+1}}{(j+1)!} (-1)^j\int_K\bb{N}^\trans \bs{S}^j\bb{N}\dd K \bs{U}_{t_i}
\end{equation*}
summing weighted spatial derivatives from zeroth to $(k+1)$-th order with different prefactors. Higher order derivatives are computed by iterative calculation of a first derivative and a subsequent projection applying the inverse mass matrix.
The Taylor--Cauchy--Kowalevski term is completely element local and significantly more calculations are operated on the read data compared to the evaluation of the inverse mass matrix. 

We propose a successive evaluation of the gradient of a field and projection instead of direct evaluation of high derivatives, which is significantly more efficient in the context of a matrix-free implementation on general meshes. Successively, problems of the form $(w,\widetilde{u}^{i+1})_K=(w,\widetilde{u}^{i})_K$ are solved, where the $j$-th spatial derivative is denoted $\widetilde{u}^j=\nabla^j u$. Algorithm~\ref{alg:naive} shows how the spatial derivatives are calculated by evaluation of the gradient field and projection onto the degrees of freedom for the non-collocated basis GL. The application of the inverse mass matrix must be understood in the matrix-free sum-factorization context according to \cite{ksmw15}.

\begin{algorithm}
\caption{Evaluation of $j$-th order derivatives for the non-collocated basis GL on general non-Cartesian grids.}
\label{alg:naive}
\begin{algorithmic}
\FOR{$i=1,\ldots,j$}
\STATE evaluate $\nabla {\widetilde u}^i_\text{GL}$ in the integration points $\bs{\xi}_\text{G}$
\STATE multiply with weighting functions evaluated in integration points $\bs{\xi}_\text{G}$ and sum to get right hand side vector $\bs{r}_\text{GL} = (w_\text{GL},\nabla \widetilde{u}^i_\text{GL})_K$
\STATE apply inverse mass matrix $\bb{M}^{-1}_\text{GL} \bs{r}_\text{GL}$ to get coefficients $\widetilde{\bs{U}}^{i+1}_\text{GL}$ of the field $\widetilde{u}^{i+1}_\text{GL}$
\ENDFOR
\end{algorithmic}
\end{algorithm}

The algorithm simplifies drastically if nodes and integration points are collocated because the weighting with test functions and the application of the inverse mass matrix cancel out and all the interpolation matrices are identity matrices. The simplified procedure is shown in Algorithm~\ref{alg:nice}. 

\begin{algorithm}
\caption{Evaluation of $j$-th order derivatives with the collocated basis G.}
\label{alg:nice}
\begin{algorithmic}
\FOR{$i=1,\ldots\,j$}
\STATE evaluate $\nabla \widetilde{u}^i_\text{G}$ in the integration points $\bs{\xi}_\text{G}$
\STATE set values for $\widetilde{\bs{U}}^{i+1}_\text{G}$ of the field $\widetilde{u}^{i+1}_\text{G}$
\ENDFOR
\end{algorithmic}
\end{algorithm}

Note that the basis change is done on the fly when processing the data from one element, 
not on the global solution vector as that would incur additional memory transfer. If the current global solution vector contains the degree of freedom values in the GL basis description, the first evaluation in Algorithm~\ref{alg:nice} involves the interpolation from GL to G and then all remaining evaluation from $i=2$ to $i=j$ are computed purely in the collocated basis G. 

\subsection{Degree Reduction}\label{sec:basis-reduction}
A possibility to further reduce the work in the
evaluation of the Taylor--Cauchy--Kowalevski sum as
presented in the previous section is to reduce the
polynomial degree in the representation of higher order spatial derivatives. This is a well-established method in the ADER-DG community for affine element shapes where typically hierarchical bases are used \cite{bhrb14,dk06}. The higher order spatial derivatives naturally give a contribution only to the lower degree coefficients of the hierarchical basis in the Taylor--Cauchy--Kowalevski procedure. In the case of sum factorization evaluation on general non-Cartesian meshes, however, the embedding to lower polynomial degrees involves additional operations as compared to the spectral evaluation routine from Algorithm~\ref{alg:nice}. In a hierarchical basis, i.e.,~Legendre polynomials on quadrilaterals or hexahedra, the integrals with a non-affine element geometry are evaluated by quadrature with sum factorization, which in turn must transform between the Legendre basis and the collocation basis. For one spatial derivative $i$ in Algorithm \ref{alg:nice}, the complexity rises from $d$ tensor product kernels per component in the spectral evaluation to $3d$ tensor product kernels, where $d$ kernels each are needed for the basis change in interpolation and integration, respectively.

A more efficient algorithm can be devised as follows. First, we reduce the degree in terms of the Lagrange basis in the respective Gaussian integration points of degree $k-i+1$ and $k-i$ of step $i$ in the Taylor--Cauchy--Kowalevski sum. The degree reduction is performed by an operation $\bb{P}_{k-i+1}^{k-i} \widetilde{u}^{i+1}$ where $\bb{P}_{k-i+1}^{k-i}$ is the projection operator from degree $k-i+1$ to $k-i$ on the reference element. Like in the hierarchical case, this setup combined with the additional interpolation of the result into the points of the Taylor--Cauchy--Kowalevski sum involves $3d$ tensor product kernels per component, as compared to only $d$ kernels for the spectral derivative. Thus, as a second ingredient we propose to apply the
degree reduction only for every second spatial derivative. This step also has the advantage of limiting the extra amount of geometry information, i.e., the inverse Jacobian, that needs to be loaded in each quadrature point, as Gauss formulas of different degrees evaluate the integrands in different positions and our implementation uses pre-computed inverse Jacobians. In other words, we re-use the data of the inverse Jacobians loaded into caches once again.
The $i$-th spatial derivative $\tilde{u}^i$ is thus expressed in
a basis of polynomial degree $k-\lfloor \nicefrac{j}{2} \rfloor \cdot 2$.
Algorithm~\ref{alg:reduc} details this procedure.
\begin{algorithm}
\caption{Evaluation of high derivatives with collocated basis $G$ and a degree reduction in every second step.}
\label{alg:reduc}
\begin{algorithmic}
\STATE set $k^{(1)}=k+1$
\FOR{$i=1,\ldots,j$}
\STATE evaluate $\nabla \widetilde{u}^i_G$ in $\left(k^{(i)}\right)^d$ integration points $\bs{\xi}^{(i)}_G$
\STATE set values for $\widetilde{\bs{U}}^{i+1}_G$ of the field $\widetilde{u}^{i+1}_G$ for degree $k^{(i)}$
\IF{i \text{mod} 2 = 0}
\STATE project $\widetilde{u}^{i+1}_G$ to degree $k-i$ by sum-factorized multiplication, $\bb{P}_{k-i+2}^{k-i}\widetilde{\bs{U}}^{i+1}_G$
\STATE set $k^{(i+1)} = k^{(i)}-2$
\ELSE
\STATE set $k^{(i+1)} = k^{(i)}$
\ENDIF
\ENDFOR
\end{algorithmic}
\end{algorithm}

\section{Performance Evaluation}\label{sec:results}
In the following subsections, the performance of DG with ADER and DG with explicit Runge--Kutta time integration is analyzed in terms of operation counts, computation time, throughput, and scalability.

If not specified otherwise, the computational setup as shown in Table~\ref{tab:broadwell} is used in the numerical examples.
\begin{table}
\begin{tabular}{|l|l|}
\hline
CPU & $2\times 14$ core Intel Xeon Broadwell E5-2690v4, 2.6 GHz 
\\ \hline 
memory & 8 channels DDR 4 (2400 MHz) $153\nicefrac{\text{GB}}{\text{s}}$ theoretic
\\ \hline 
compiler & \texttt{g++} version 6.2 
\\ \hline 
compiler optimization & \texttt{-march=haswell -03 -funroll-loops}
\\ \hline
\end{tabular}
\caption{System specifications for the numerical tests.}
\label{tab:broadwell}
\end{table}

\subsection{Operation Counts}
For the ADER-DG method as given in equation~\eqref{eq:ader}, there are three main contributions to the computational costs, namely the application of the inverse mass matrix cost $C_\bb{M}$, the application of the global derivative operator cost $C_\bb{K}$, and the evaluation of the Taylor--Cauchy--Kowalevski sum cost $C_\text{TCK}$, resulting in an overall cost of
\begin{align}
C_\text{ADER-DG} = 2 \cdot C_\bb{M} + C_\bb{K} + C_\text{TCK}.
\end{align}
For an $s$-stage Runge--Kutta scheme, the costs are
\begin{align}
C_\text{RK} = s\cdot (C_\bb{M} + C_\bb{K}).
\end{align}
In ADER, the high-order approximations contribute to a sum over all terms in the truncated Taylor series and $C_\text{TCK}$ involves an additional dependency on the polynomial degree $k$ compared to $C_\bb{M}$ and $C_\bb{K}$. In contrast, high-order approximations with Runge--Kutta schemes use more stages $s$, where the number of stages $s$ has to increase faster than the required order of accuracy (for orders larger than four) due to the Butcher barriers.
ADER-DG as well as Runge--Kutta methods both repeatedly call the application of the derivative operator and of the inverse mass matrices. The main difference is that ADER-DG applies the derivative operator locally, i.e., element-wise, while Runge--Kutta integrators rely on the global derivative operator containing both element and face contributions.

The operation counts for $C_\bb{M},C_\bb{K},$ and $C_\text{TCK}$ are derived from vector updates, matrix-vector or matrix-matrix multiplications, which in turn rely on the matrix-free implementation of integral evaluation for tensorial shape functions explained in Section~\ref{sec:algo}. The cost for the evaluation of one tensor product kernel on a $\delta$-dimensional domain (e.g. $\delta=d$ in an element or $\delta = d-1$ on element faces) calculates to
\begin{align*}
C_\text{tensorial}(\delta) = \left( 2 \cdot \frac{k+1}{2} \cdot 2 + (k+1) + 2 \cdot  \left\lfloor\frac{(k-1)\cdot(k+1)}{2} \right\rfloor\right) \cdot (k+1)^{\delta-1},
\end{align*}
where the three summands in braces represent additions and subtractions (first term), multiplications (second term), and fused multiply-add operations (last term, counted as two arithmetic instructions). Note that an even-odd decomposition of the local coefficients and matrices is used that cuts operation count into approximately one half compared to a usual 1D matrix-vector product \cite{kopriva,kk17}.
This cost is the main building block to derive the operation counts for $C_\bb{M},C_\bb{K},$ and $C_\text{TCK}$ with the number of calls to the tensor product kernels as summarized in Table~\ref{tab:kernel-calls}. For details on the derivation of operation counts, we refer to~\cite{kk17} and~\cite{ksmw15}.

\begin{table}
\centering
\caption{Number of calls of the tensor product kernel for acoustics in terms of cell kernels "C" and face kernels "F".}
\label{tab:kernel-calls}
\begin{tabular}{|l|lll|}
\hline
& mass $\bb{M}^{-1}$ & stiffness matrix $\bb{K}$ & TCK \\ \hline
cell eval. & $2d$ &  $d^2+2d$ & $2d$ \\
cell deriv. & $-$ & $2d$ & $(k-1)\cdot d$  \\
face eval. & $-$ & $d\cdot 4(d-1)$ & $-$ \\ \hline
total & $2d\cdot \text{C}$ &  $(d^2+4d)\cdot\text{C}+4(d^2-d)\cdot\text{F}$ & $(2d+d(k-1))\cdot\text{C}$ \\ \hline
\end{tabular}
\end{table}

\begin{figure}
\centering
\includegraphics[scale=1.0]{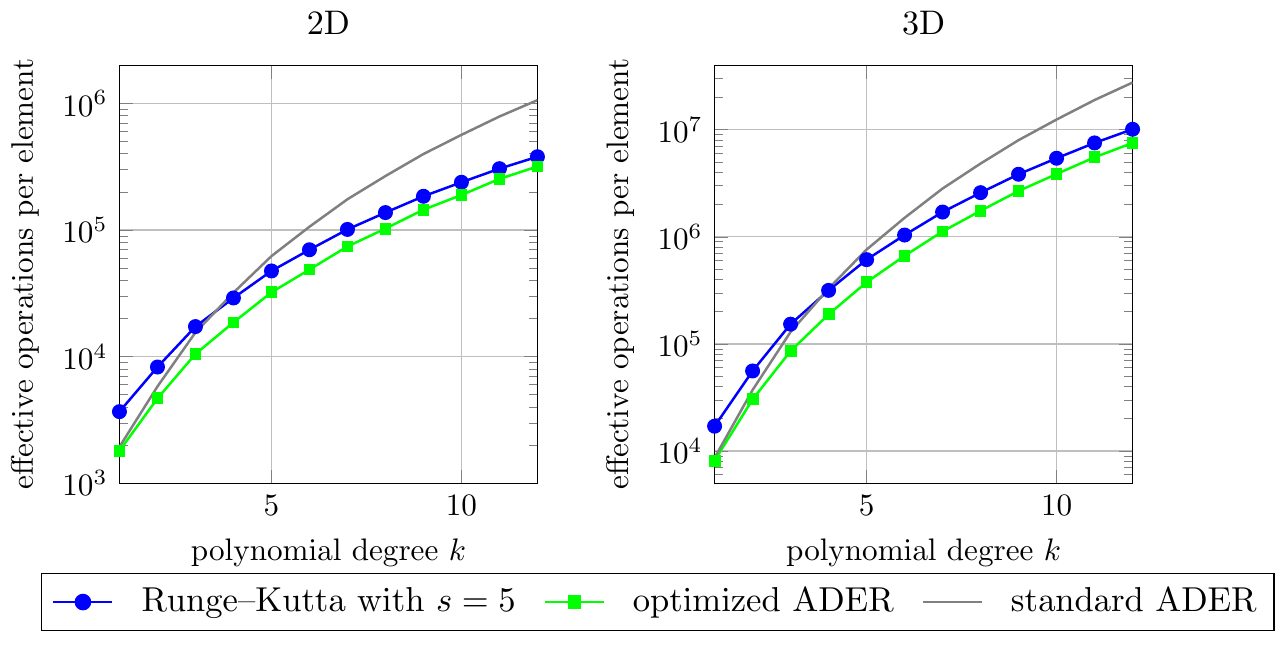}
\caption{Operation counts for evaluation of one element with a five stage Runge--Kutta scheme and ADER-DG  in two and three dimensions.}
\label{fig:comp-opcounts-RK-ADER}
\end{figure}

In Fig.~\ref{fig:comp-opcounts-RK-ADER}, operation counts for the full schemes are compared, i.e., ADER versus a Runge--Kutta scheme with five stages for all polynomial degrees. To allow for generalization and allow comparability, the number of stages for the Runge--Kutta scheme is kept constant. 
ADER involves fewer arithmetic operations for all considered polynomial degrees $k\in\{1,\ldots,12\}$ for both $d=2,3$. The figure also shows operation counts in case no basis change to a collocated basis $G$ and no degree reduction for the higher order spatial derivatives are carried out, i.e., if the proposition of Sections~\ref{sec:basis-change} and~\ref{sec:basis-reduction} are not realized in the implementation of the Taylor--Cauchy--Kowalevski procedure. Without optimizations, it is apparent that ADER has a higher polynomial dependency on $k$ than Runge--Kutta. The optimizations however compensate this dependency. The basis change is not applied to the Runge--Kutta discretization because the element-local mass matrix inversion is dominated by the memory bandwidth in contrast to the ADER Taylor--Cauchy--Kowalevski term with higher operational cost, see also Fig.~\ref{fig:rk4-basis}.

\begin{figure}
\centering
\includegraphics[scale=1.0]{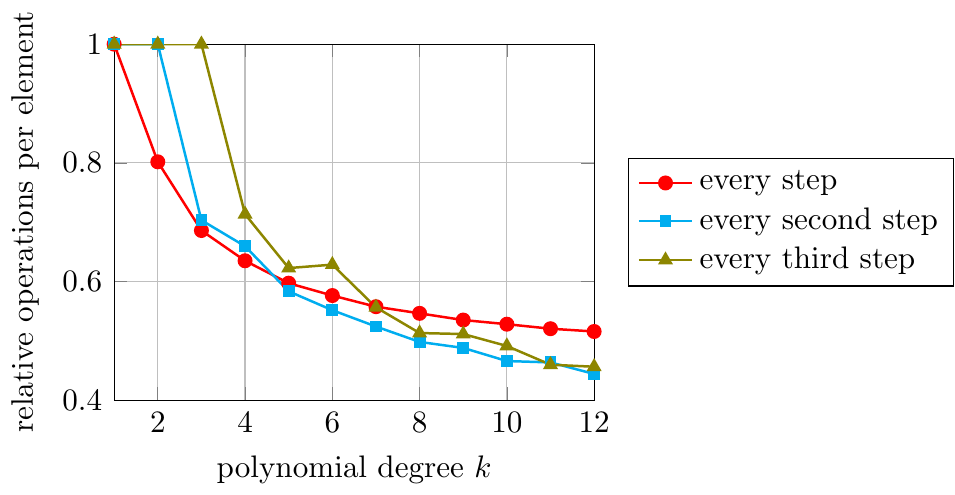}
\caption{Operation counts for the Taylor--Cauchy--Kowalevski procedure with degree reduction every step, every second, and every third step relative to using the degree $k$ for all summands.}
\label{fig:reduc-comp}
\end{figure}

Fig.~\ref{fig:reduc-comp} visualizes the gains of the degree reduction approach. It plots the operation counts for the Taylor--Cauchy--Kowalevski procedure with degree reduction relative to an implementation without degree reduction and compares the reduction in every step, every second step, and every third step. The operational cost is halved. For higher orders, the reduction in every step is not preferable because the projection introduces overhead.  For moderate polynomial degrees, the approach to reduce the basis in every second step of the Taylor--Cauchy--Kowalevski procedure appears most beneficial.

\subsection{Computational Timings}
A two and a three dimensional example are set up on a domain $\Omega=[0,1]^d$, solving the acoustic wave equation with vibrational modes as an analytic solution as in~\cite{skw17}. We test on a mesh with slightly deformed elements to prevent the built-in Cartesian mesh optimizations in the code of \cite{kk17} for polynomial degrees $k\in\{1,\ldots,12\}$. In 2D, $1280^2$ and $640^2$ elements are used for $k\in\{1,2,3\}$ and $k\in\{4,\ldots,12\}$, respectively.
In 3D, $80^3$ and $40^3$ elements are used for the respective polynomial degrees.

\begin{figure}
\centering
\includegraphics[scale=1.0]{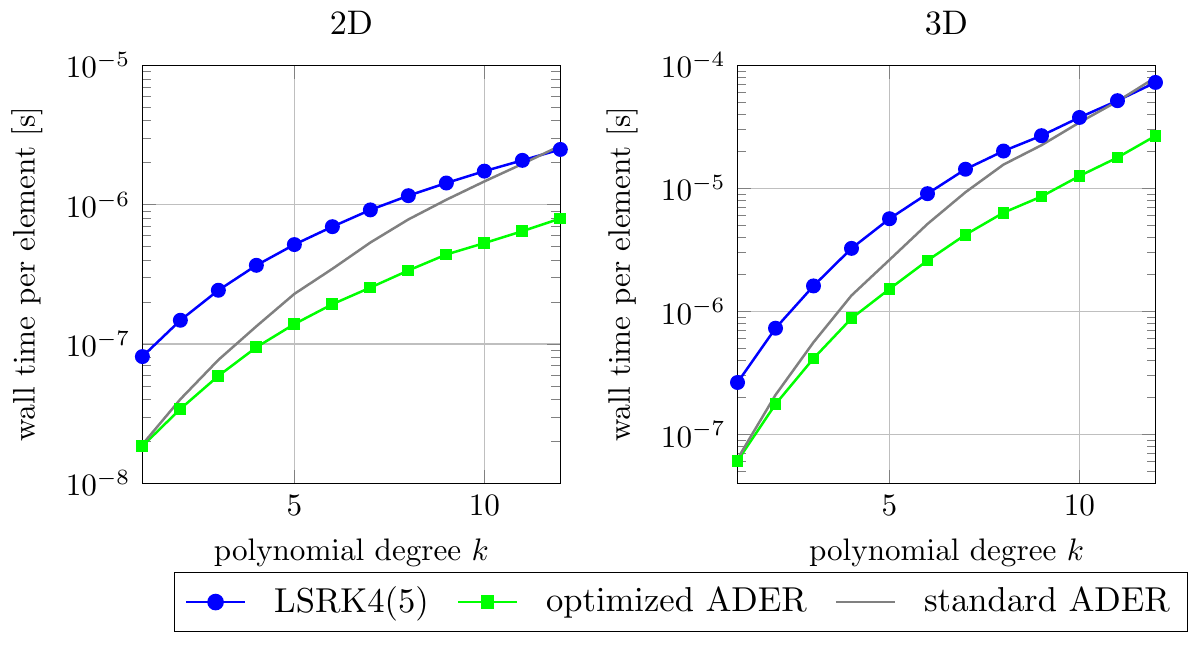}
\caption{Run time per element for an explicit Runge--Kutta scheme with $s=5$ and ADER-DG in two and three dimensions, measured on 28 cores.}
\label{fig:comp-time-RK-ADER}
\end{figure}

Fig.~\ref{fig:comp-time-RK-ADER} shows measured run times per time step and per element for numerical experiments on $28$ cores with LSRK4(5). 
Comparison between Fig.~\ref{fig:comp-opcounts-RK-ADER} and  Fig.~\ref{fig:comp-time-RK-ADER} reveals that the run time follows the operation counts.
ADER is significantly more efficient in terms of wall time per element, though, which is due to the fact that the  Taylor--Cauchy--Kowalevski procedure requires one global vector access but Runge--Kutta requires one global vector access for each application of the  global derivative operator $\bb{K}$. 
In other words, the operation counts that ignore the memory access fail to accurately predict the run time. This statement also holds for the ADER implementation without basis change with run times between LSRK4(5) and the optimized ADER.

\subsection{Breakdown into Algorithmic Components}\label{sec:alg-comps}
\begin{figure}
\centering
\includegraphics[scale=1.0]{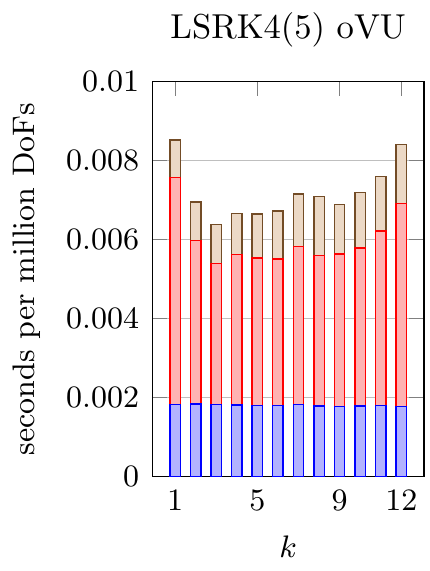}
\includegraphics[scale=1.0]{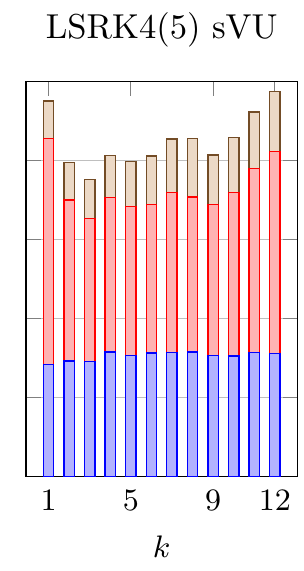}
\includegraphics[scale=1.0]{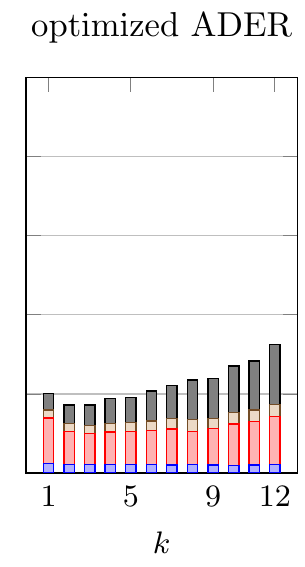}
\includegraphics[scale=1.0]{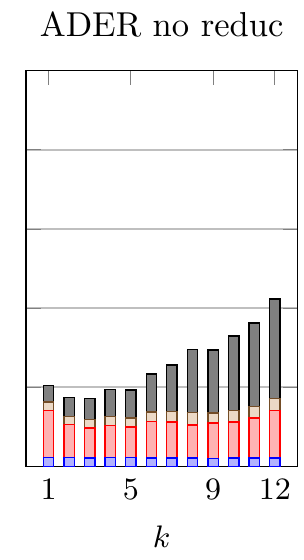}
\includegraphics[scale=1.0]{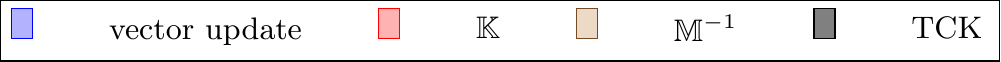}
  \caption{The algorithmic components for various polynomial degrees in 3D on 28 cores with the low-storage scheme LSRK4(5) with optimized vector update ``oVU'' routine and standard vector update ``sVU'', and ADER time integration with the two optimizations as presented in Sections~\ref{sec:basis-change} and~\ref{sec:basis-reduction}. In the very right panel, results are shown for the ADER algorithm without the degree reduction optimization.}
  \label{fig:alg-comps-3d}
\end{figure}

Let us examine the algorithmic components separately in terms of the computational time per million degrees of freedom. 
In Fig.~\ref{fig:alg-comps-3d}, we compare two versions of LSRK4(5) with an optimized vector updater that only reads and writes two vectors per stage by merging the loop over vectors into a single loop and a standard vector updater reading five vectors and writing two vectors per stage, similar to the ``merged operations'' presented in~\cite{eg16}. Also, the results for the optimized ADER scheme implementing the basis change and degree reduction from Sections~\ref{sec:basis-change} and~\ref{sec:basis-reduction} and the ADER scheme without the degree reduction are shown in Fig.~\ref{fig:alg-comps-3d}.
The results are obtained from simulations on a three-dimensional non-Cartesian grid. The measured run time is the accumulated time spent in the respective functions, e.g., the contribution of the derivative operator $\bb{K}$ appears five times as high for LSRK4(5) compared to ADER, because it is applied in each of the five stages while ADER applies the derivative operator only once per time step.

As can be seen from Fig.~\ref{fig:alg-comps-3d}, the cost for vector updates in Runge--Kutta schemes cannot be neglected, neither in the standard implementation nor in the optimized variant, where they contribute with about a third or a quarter of the run time, respectively. Likewise, our results document the high level of performance reached in the computations of $\bb{K}$ and $\bb{M}^{-1}$, a distinctive feature of our implementation. For ADER, the vector updates have a comparably small contribution to the entire cost because only one single update at the end of the method is required.

The application of the derivative operator $\bb{K}$ is most expensive for $k=1$ with its rather disadvantageous ratio between degrees of freedom located on the element faces as compared to the interior. The ratio improves for higher orders and an almost constant throughput per degree of freedom is obtained, despite the theoretical $\mathcal O(k)$ increase in arithmetic complexity. The constant throughput is mainly explained by the memory transfer that scales as $\mathcal O(1)$ per unknown. The Taylor--Cauchy--Kowalevski procedure shows slightly increasing costs for higher degrees as an additional summand contributes in the Taylor expansion according  to equation~\eqref{eq:ader}. Nonetheless, the increase for higher degrees is moderate due to the proposed degree reduction approach as presented in Section~\ref{sec:basis-reduction} and efficient algorithms, less than doubling the run time between degree two and twelve. In the rightmost panel of Fig.~\ref{fig:alg-comps-3d}, results are shown for a   Taylor--Cauchy--Kowalevski procedure that does not reduce the polynomial degree, where the increase in computing time is much more significant. Obviously, the latter reaches higher arithmetic throughput with more than 300 GFLOPs/s, which is a secondary quantity, though.

\subsection{Roofline Performance Model}
\begin{figure}
\centering
\includegraphics[scale=1.0]{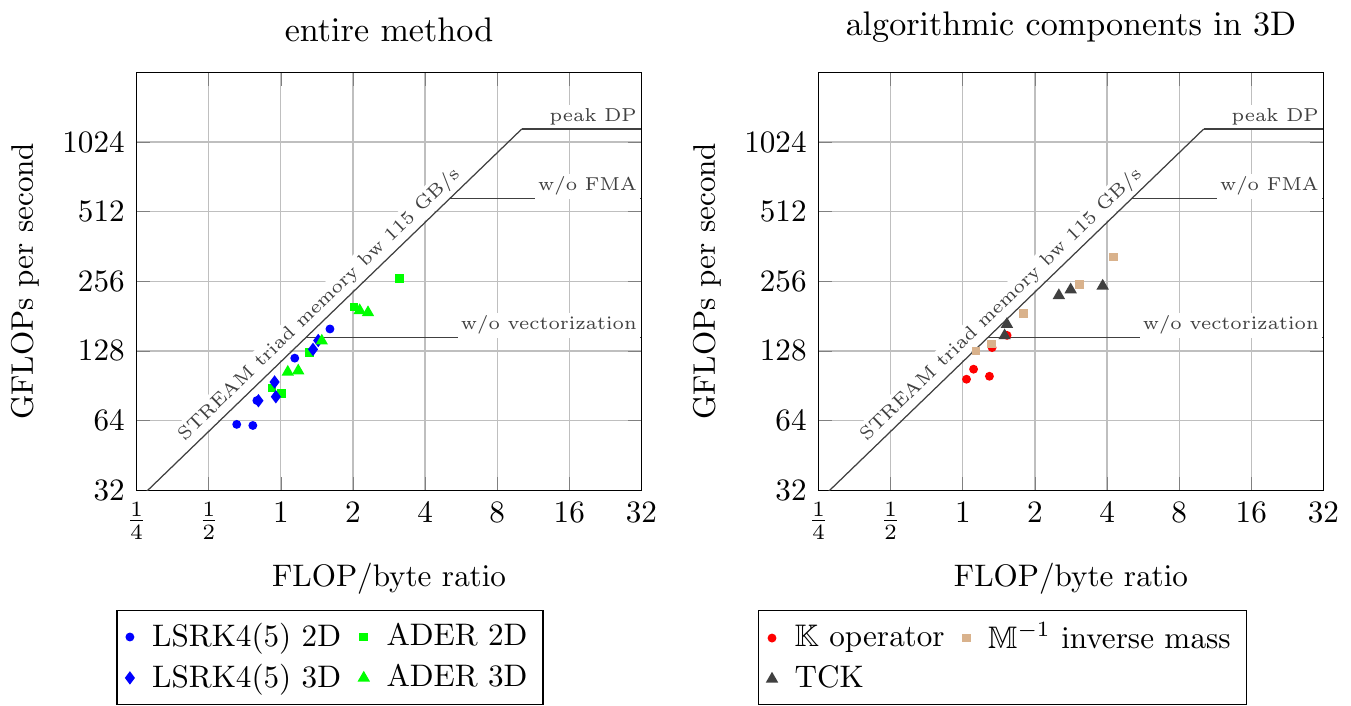}
\caption{Roofline model for polynomial degrees $k=1,2,4,8,12$ for the entire method (left panel) and the individual method components in 3D (right panel).}
\label{fig:roofline}
\end{figure}

Fig.~\ref{fig:roofline} shows a roofline model according to~\cite{roofline} for LSRK4(5) and ADER in two and three dimensions on a non-Cartesian mesh. The roofline plots contain the hardware specific limits in terms of the memory bandwidth limit (diagonal lines) measured with the STREAM triad benchmark and the peak arithmetic performance (horizontal lines). All numbers are based on measured data from hardware performance counters, extracted from monitoring our programs with the \texttt{likwid} performance measurement tool, version 4.3, as presented in~\cite{psti}.

Generally, the polynomial degree increases for the points from left to right. The left panel of Fig.~\ref{fig:roofline} highlights that ADER comes with a higher arithmetic intensity, in particular for the higher degrees, as compared to LSRK4(5) that is clearly in the memory bandwidth bound regime.  The computations are made for a non-Cartesian mesh where not only vector entries and some index data must be loaded from main memory but also the inverse of the Jacobian transformation in each quadrature point.

The right panel of Fig.~\ref{fig:roofline} shows the results for the individual components of the methods. The evaluation of the derivative operator $\bb{K}$ that needs to load geometry data is most strongly limited by the memory. Note that loops over cells and faces are interleaved in our implementation to re-use the vector data loaded into caches in cell integrals also for face integrals. More detailed measurements similar to the ones presented in \cite{kk17} show that the code of face integrals finds more than 90\% of the vector data for face integrals already in caches. Nonetheless, the partial indirect addressing with gather/scatter type instructions to rearrange the face data for vectorization over several faces and the remaining cache misses reduce throughput by around 25\% as compared to idealized code that only performs the integration, see the experiments in \cite{kk17} for details. The Taylor--Cauchy--Kowalevski procedure comes along with considerably higher arithmetic intensities, but the high level of arithmetic optimizations in the proposed algorithm keeps the code still in the memory-limited region, in clear contrast to e.g.~\cite{bhrb14}. Given these detailed results, the behavior of the entire method can be characterized as follows: the Runge--Kutta scheme applies $\bb{K}$ and $\bb{M}^{-1}$ five times in each time step, while ADER only applies $\bb{K}$ and $\bb{M}^{-1}$ once and the rest is taken care of by the completely element-local Taylor--Cauchy--Kowalevski procedure, enabling a much higher re-use of cached data.

\subsection{Throughput as Function of Polynomial Degree}
\begin{figure}
\centering
\includegraphics[scale=1.0]{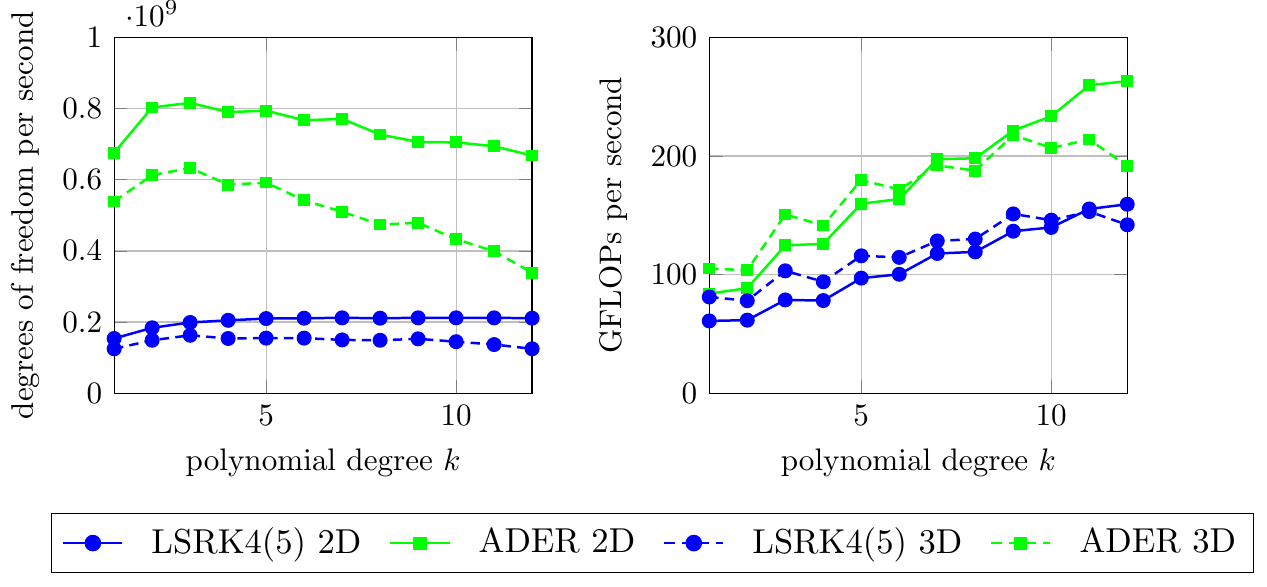}
\caption{Throughput in terms of degrees of freedom per second and GFLOPs per second comparing LSRK4(5) and ADER-DG in two and three dimensions.}
\label{fig:comp-time-throughput}
\end{figure}

Fig.~\ref{fig:comp-time-throughput} lists the computational throughput in terms of degrees of freedom processed per second and the actually realized GFLOPs per second rate, measured with the \texttt{likwid} tool \cite{psti}, version 4.3. If not stated otherwise, the optimized implementation of ADER using the basis change and reduction as proposed in Sections~\ref{sec:basis-change} and~\ref{sec:basis-reduction} are considered.

The performance advantage of ADER compared to the low storage Runge--Kutta scheme is a factor of around $4$ for low polynomial degrees. The advantage decreases to a factor of $1.5$ for $k=12$ owing to the additional computations for the high order of accuracy of the ADER time integration, whereas the LSRK4(5) schemes uses a fixed temporal accuracy of four with 5 stages. In terms of GFLOPs rates, which lie between $60$ and $264$ GFLOPs per second, ADER operates $1.6$ times more GFLOPs per second due to better caching. For comparison, a matrix based evaluation of a sparse matrix vector product was reported to reach $21$ GFLOPs per second~\cite{Kronbichler16b} on the same hardware. 

\subsection{Throughput as Function of Problem Size}
\begin{figure}
\centering
\includegraphics[scale=1.0]{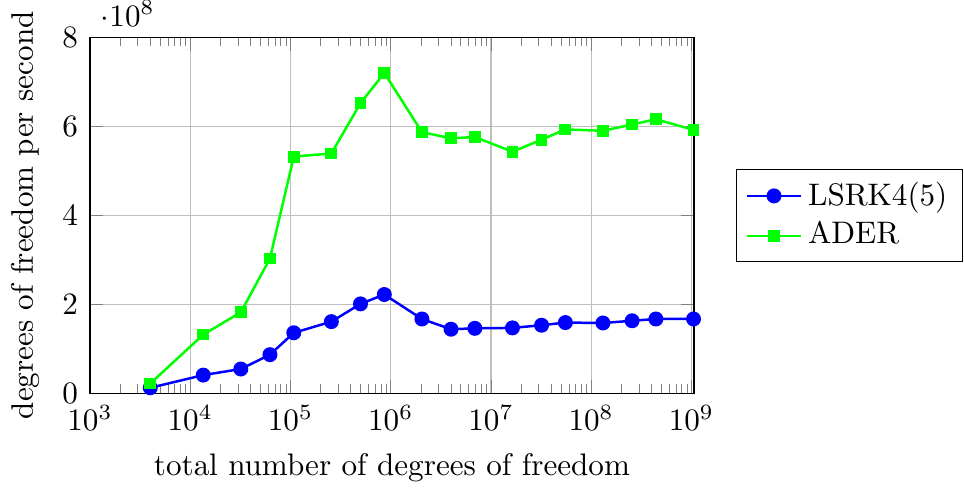}
\caption{Throughput as function of the problem size, displayed as the number of degrees of freedom processed per second for one time step for $k=4$.}
\label{fig:throughput-problemsize}
\end{figure}

The throughput as a function of the problem size between $4\cdot10^3$ and $2\cdot10^9$ degrees of freedom in 3D with shape functions of polynomial degree $k=4$ is shown in Fig.~\ref{fig:throughput-problemsize}. For small discretizations of size $n^\text{dof}<10^5$, limited parallelism prevents the full exploitation of the 28 cores. Then, the parallel efficiency improves and a first peak is reached at around $n^\text{dof}\approx 10^6$ where all data fits into the 70 MB of level 3 cache of the two processors and no access to main memory is needed. Performance drops again once the data structures exceed the caches and the vector and geometry data need to be streamed from main memory. For $n^\text{dof}>10^6$, a slow increase of the throughput is noted, which is related to the decreased influence of the MPI communication for larger discretizations.

\subsection{CPU Time Versus Accuracy}\label{sec:time-vs-acc}

\begin{figure}
\centering
\includegraphics[scale=1.0]{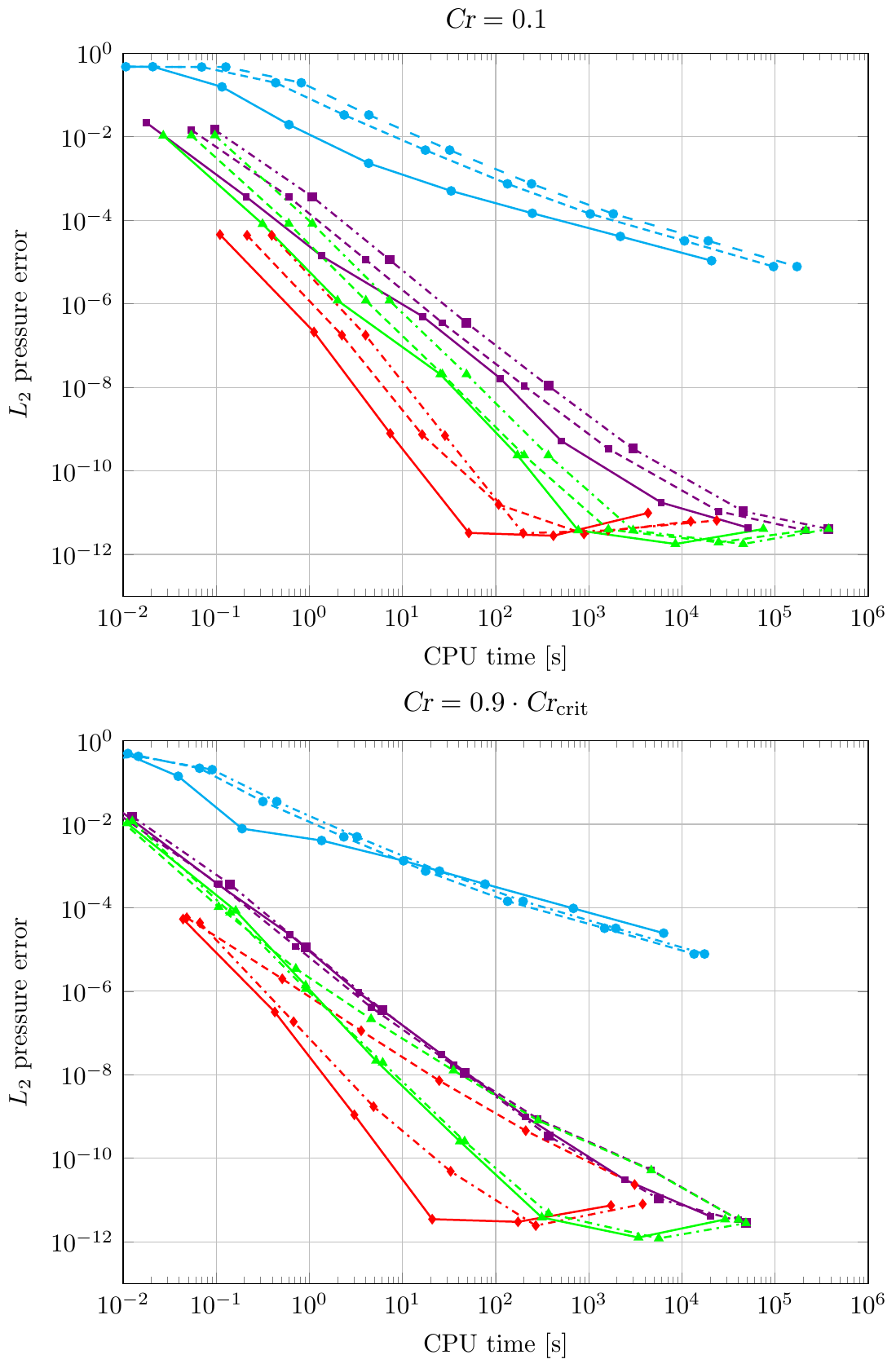}
\includegraphics[scale=1.0]{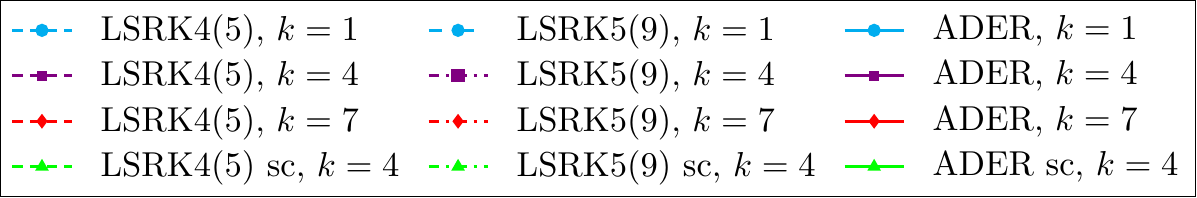}
\caption{Accuracy over wall time accumulated over all MPI ranks comparing low storage Runge--Kutta and ADER for the same time step size and at their respective critical time step size. Polynomial degrees $k=1,4,7$ are studied. For $k=4$ the superconvergent pressure solution is also considered for error calculation indicated by ``sc''.}
\label{fig:solv-acc-RK-vs-ADER-post}
\end{figure}

In this section, the time to solution is evaluated with respect to the accuracy on a two-dimensional example. With the standard setup of a vibrating membrane with seven modes, simulations are run on different refinement levels of the mesh between $5^2$ and $1280^2$ elements for polynomial degrees $k=1,4,7$. 
For a fair comparison, the number of processors is chosen reasonably for the discretization sizes, e.g., the coarsest discretization is computed on one processor core only while the finest discretization is computed on $28$ cores. The final numbers report the accumulated CPU time over all utilized processors. For the accuracy, the $L_2$ pressure error at time $t^\text{end}=1.0$ is considered. The upper panel of Fig.~\ref{fig:solv-acc-RK-vs-ADER-post} plots the results for a Courant number of $\textsl{Cr} = 0.1$. Generally, higher orders appear beneficial in case a strict accuracy criterion is applied. ADER yields the same solution quality as LSRK4(5) in less computational time.

Since ADER and Runge--Kutta schemes are subject to different CFL stability limits, we evaluate the critical Courant number $\textsl{Cr}_\text{crit}$ by an iterative procedure. The results are summarized in Table~\ref{tab:crit-Cr} for the considered two-dimensional as well as for a three-dimensional setup.
\begin{table}
\centering
\caption{Critical Courant numbers $\textsl{Cr}_\text{crit}$ for LSRK4(5), LSRK5(9) and ADER determined by an iterative procedure on a two- and three-dimensional domain with Cartesian mesh.}
\label{tab:crit-Cr}
\begin{tabular}{c|ccc|ccc}
& \multicolumn{3}{c|}{2D} & \multicolumn{3}{c}{3D} \\ \hline
 & LSRK4(5) & LSRK5(9) & ADER  & LSRK4(5) & LSRK5(9) & ADER\\ \hline
$k=1$ & $0.44$ & $0.55$ & $0.18$ & $0.28$ & $0.36$ & $0.13$ \\
$k=4$ & $0.69$ & $0.87$ & $0.28$ & $0.49$ & $0.58$ & $0.19$ \\
$k=7$ & $0.68$ & $0.88$ & $0.28$ & $0.44$ & $0.57$ & $0.18$
\end{tabular}
\end{table}
The critical Courant number for ADER is approximately $2.4$ and $3.1$ times smaller compared to the low storage Runge--Kutta schemes LSRK4(5) and LSRK5(9) from \cite{kcl00}. Therefore, we also analyze the solution accuracy versus the CPU time at $\textsl{Cr} = 0.9\cdot \textsl{Cr}_\text{crit}$ of the respective time integrator in the bottom panel of Fig.~\ref{fig:solv-acc-RK-vs-ADER-post}. For $k=1$, ADER is faster for most error tolerances. For $k=4$, all methods perform similarly. With polynomial degree $k=7$ for the spatial discretization, ADER outperforms LSRK4(5) clearly because the LSRK4(5) scheme is of order four while the spatial discretization is of order eight, and the temporal error dominates over the spatial error. This is also true for LSRK5(9) at high spatial resolution.

In Section~\ref{sec:hdg}, HDG was introduced as a special DG method with a superconvergence property. By a simple postprocessing step, a superconvergent pressure solution is constructed that converges with rate $k+2$. The results for $\textsl{Cr}=0.1$ and $\textsl{Cr} = 0.9\cdot \textsl{Cr}_\text{crit}$ for polynomial shape functions of degree $k=4$ are presented as green lines in Fig.~\ref{fig:solv-acc-RK-vs-ADER-post}. Again, for $\textsl{Cr}=0.1$, ADER performs slightly better compared to LSRK4(5), while a significant performance advantage for ADER is noted for $\textsl{Cr} = 0.9\cdot \textsl{Cr}_\text{crit}$, which is due to better temporal accuracy. Comparing the top and bottom panel for $k=7$ shows that ADER is also faster than LSRK4(5) with smaller time steps.  The discretization with LSRK5(9) is competitive to ADER for the superconvergent pressure result with $k=4$. For the polynomial degree $k=7$, a change in slope indicates the turnover from spatial error domination to temporal error domination. Where the temporal error dominates, the advantage for ADER is more distinct.
This is due to the fact that ADER-DG is automatically $k+1$ convergent in space and time while Runge--Kutta is limited by the Butcher barriers: either overproportionally more stages are required to match the temporal order of accuracy with the spatial discretization, or a small time step is required.

\subsection{Scalability}
\begin{figure}
\centering

\includegraphics[scale=1.0]{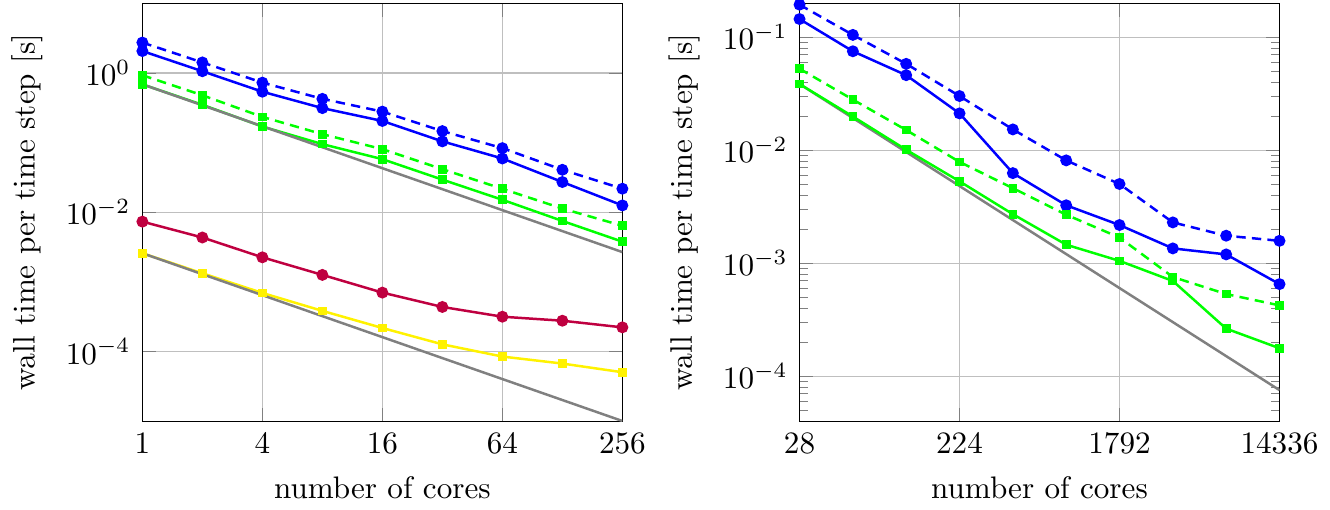}
\includegraphics[scale=1.0]{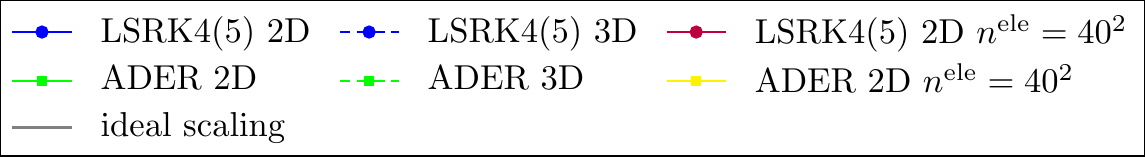}
\caption{Strong scaling for LSRK4(5) and ADER-DG in two and three dimensions in terms of wall time per time step over the number of cores. The left panel is obtained on a cluster using $1$ to $256$ cores of Intel Haswell E5-2630 v3, while the right panel is obtained with simulations on the SuperMUC Phase 2 system using $28$ to $14336$ cores of Intel Haswell E5-2697 v3.}
\label{fig:scaling}
\end{figure}

In order to assess the strong scalability of the proposed methods, a two and a three dimensional geometry consisting of $640^2$ and $40^3$ elements of polynomial degree $k=4$ ($n^\text{dof}=3.1\cdot 10^7$ and $n^\text{dof}=3.2\cdot 10^7$, respectively) are used and the solution is computed on $1$ to $2^8=256$ processors on a parallel cluster of $2\times 8$ core Intel Xeon E5-2630 v3 (Haswell) processors at 2.4 GHz. Additionally, one computation with fewer elements ($40^2$, $n^\text{dof}=1.2\cdot 10^5$) in 2D is carried out such that the distribution yields only six elements per processor for the highest level of parallelism. The left panel of Fig.~\ref{fig:scaling} summarizes the results.  In accordance with the previous sections, ADER is consistently faster than LSRK4(5) for the same time step size $\Delta t$. The scaling is almost ideal with a slight kink when going from $8$ to $16$ processors, where the code goes from being compute bound to being memory bound. For the small discretization the scaling deteriorates for high processor numbers due to communication overhead. In the right panel, results are shown for a strong scaling study on the SuperMUC Phase 2 Petascale system with nodes of $2\times 14$ Intel Xeon E5-2697 v3 (Haswell) processors at 2.6~GHz on $1$ to $512$ nodes. For the 2D LSRK4(5) simulation, a kink can be seen when going from $224$ to $448$ cores indicating that simulations are memory bandwidth bound on lower core counts but get computation bound at higher core counts. The scaling is close to ideal considering that the highest level of parallelism corresponds to $28$ and $4$ elements per processor in 2D and 3D, respectively.

\begin{figure}
\centering

\includegraphics[scale=1.0]{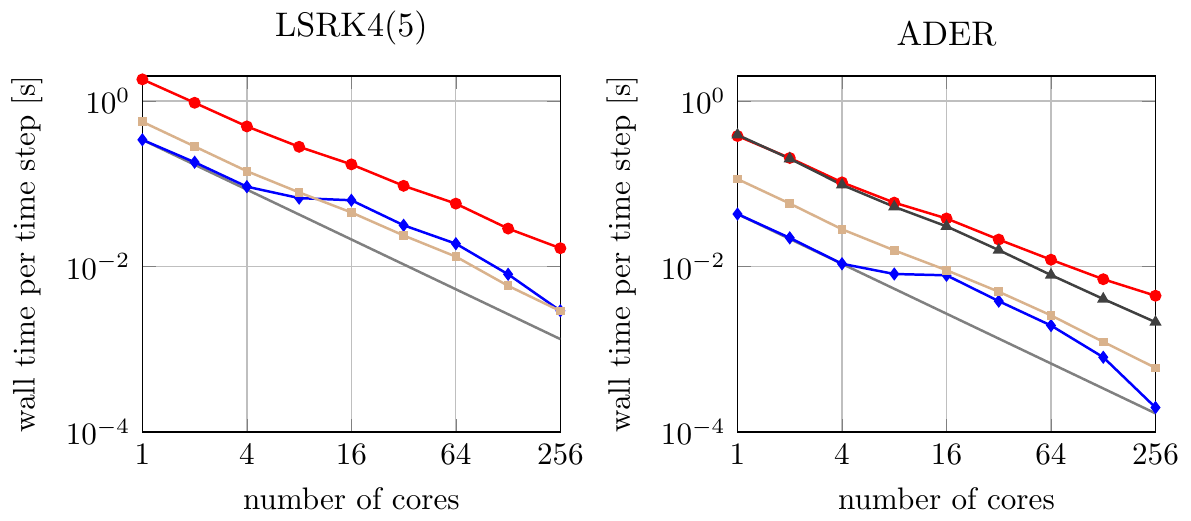}
\includegraphics[scale=1.0]{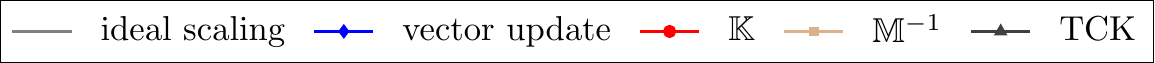}
\caption{Strong scaling for LSRK4(5) and ADER-DG in three dimensions for the components of the algorithm in terms of wall time per time step over the number of cores for $n^\text{dof} = 3.2\cdot 10^7$ on up to 16 Haswell E5-2630 v3 nodes.}
\label{fig:scaling3d-comps}
\end{figure}

Fig.~\ref{fig:scaling3d-comps} plots the strong scalability for the algorithmic components. It can be seen that the vector updates give the largest contribution to the reduced scalability between $8$ and $16$ processors which is due to the fact that the utilized cluster possesses $16$ cores but only $8$ memory channels which can be saturated already with $8$ processors. The kink is only due to shared memory effects because the communication is considered as part of the operator application. Since LSRK4(5) spends a larger fraction of time in vector updates, the reduced scalability in the overall method is more pronounced. The inverse mass matrix and the Taylor--Cauchy--Kowalevski procedure scale almost perfectly while the stiffness matrix shows a slight scaling decay due to a worse volume-to-surface effect of the data that must be exchanged with MPI.

In order to illustrate the effect of memory bandwidth on the algorithmic components, we repeat the calculations on an Intel Xeon Phi 7210-F (Knights Landing, KNL) system with 64 cores and 16 GB of high-bandwidth memory delivering up to 420 GB/s. For KNL, the code is vectorized with AVX-512, i.e., eight-wide SIMD lanes, as opposed to AVX2 with four-wide SIMD lanes on Broadwell. In Table~\ref{tab:comp-supermuc-knl}, the computing times per time step are compared to a calculation on $28$ Broadwell cores. The KNL system yields a speed up for the vector update of about $4$ for LSRK4(5), corresponding to the four times higher memory bandwidth. Also, the application of the inverse mass matrix is substantially faster. The derivative operator $\bb{K}$, on the other hand, runs slightly slower on KNL because of the more irregular code patterns in face integrals that favor the more sophisticated CPU cores of Broadwell.
\begin{table}
\caption{Comparison of one node with 28 Broadwell cores and 64 KNL cores in terms of wall time per time step for a system with $3.2\cdot10^7$ degrees of freedom. ``VU'' abbreviates ``vector update''. }\label{tab:comp-supermuc-knl}
\begin{tabular}{l|cccc|c}
 & VU & $\bb{K}$ & $\bb{M}^{-1}$ & TCK & sum \\ \hline
LSRK4(5) Broadwell & $5.4\cdot 10^{-2}$ & $1.1\cdot10^{-1}$ & $2.6\cdot10^{-2}$ & --- & $1.9\cdot10^{-1}$ \\ 
LSRK4(5) KNL & $1.3\cdot 10^{-2}$ & $1.2\cdot10^{-1}$ & $2.0\cdot10^{-2}$ & --- & $1.5\cdot10^{-1}$  \\
ADER Broadwell & $4.5\cdot 10^{-3}$ & $2.5\cdot10^{-2}$ & $6.6\cdot10^{-3}$ & $1.8\cdot10^{-2}$ & $5.4\cdot 10^{-2}$ \\ 
ADER KNL & $1.5\cdot 10^{-3}$ & $2.7\cdot10^{-2}$ & $3.6\cdot10^{-3}$ & $1.7\cdot10^{-2}$ & $4.9\cdot 10^{-2}$\\
\end{tabular}
\end{table}

\section{Conclusion}\label{sec:concl}

We presented a performance analysis for explicit Runge--Kutta and ADER DG implementations. ADER outperforms optimized low storage Runge--Kutta schemes over a range of test scenarios. The ingredients are fast integration techniques with sum factorization that combine optimal-complexity mathematical algorithms utilizing the tensor product structure of the shape functions with a highly competitive implementation that vectorizes over several elements and faces. The methods have been devised to be applicable also for complex meshes with curved quadrilateral or hexahedral elements. Our experiments clearly show that it is most efficient to evaluate operators including face integrals with nodal basis functions which have nodes on the element boundaries, while the cell evaluations in the Taylor--Cauchy--Kowalevski procedure of ADER are best performed with Lagrange polynomials in the points of the quadrature formula. To combine these two, an on-the-fly change between the bases has been proposed in this work. This result is in contrast to the consensus belief in spectral elements that favor collocation of the polynomials nodes and quadrature points.
While the theoretically derived operation counts already signify a slight benefit for ADER compared to Runge--Kutta, the actual timings show a distinct benefit, reducing the time to perform one time step by approximately a factor of four, because ADER suits modern hardware architecture better: while Runge--Kutta schemes are mostly limited by the memory bandwidth, ADER performs more operations on the data that is loaded from main memory and thus reaches a higher arithmetic intensity. A detailed analysis of Runge--Kutta versus ADER integration at the CFL stability limit has shown comparable performance where the Runge--Kutta time discretization order matches the spatial discretization order. For approximations with a high order of accuracy where the Butcher barriers set in, ADER  exceeds the abilities of Runge--Kutta because its computational cost does not grow overproportionally. 
While the findings for ADER are limited to linear hyperbolic PDEs, the optimizations regarding the basis functions and reduced vector access for the Runge--Kutta time integrators regarding basis functions are also directly applicable to general nonlinear systems of hyperbolic PDEs.

Our work highlights the importance to develop modern DG solvers according to the trends and limits in modern hardware architectures. For common solvers, the memory bandwidth limit is more relevant also when performing the relatively expensive computations of high order DG methods, i.e., is met earlier than the arithmetic performance limit. Trading global for element-local operations counters this effect, rendering approaches as the Taylor--Cauchy--Kowalevski procedure favorable. Vector updates are inherently memory bandwidth limited and need to be optimized specifically. Our experiments and performance models highlight that going significantly beyond the throughput recorded in this work demands either hardware with higher memory bandwidth, such as GPUs or the Xeon Phi, or new software paradigms that reduce the memory access over several stages, such as wavefront blocking that is already commonly used in the finite difference community.

\bibliographystyle{abbrv}
\bibliography{reference_jouabbriv}

\end{document}